\documentclass[12pt,oneside,english]{amsart}
\usepackage[latin1]{inputenc}
\usepackage{amssymb}

\makeatletter
\newtheorem{theorem}{Theorem}

\newtheorem{definition}{Definition}

\usepackage{babel}

\makeatother
\begin{document}
\baselineskip=17pt

\title{On Ramanujan Cubic Polynomials}

\author{Vladimir Shevelev}

\address{Departments of Mathematics \\Ben-Gurion University of the
 Negev\\Beer-Sheva 84105, Israel. e-mail:shevelev@bgu.ac.il}


\begin{abstract}
A polynomial $x^3+px^2+qx+r$ with the condition $pr^{\frac 1
3}+3r^{\frac 2 3}+q=0$ we call a Ramanujan cubic polynomial (RCP).
We study different interest properties of RCP, in particular, an
important role of a parameter $\frac{pq}{r}$. We prove some new
beautiful identities containing sums of 6 cubic radicals of values
of trigonometrical functions as well.

\slshape
Paper is deducated to the $120-th$ anniversary of Srinivasa
Ramanujan.\upshape
\end{abstract}

\maketitle
\section{INTRODUCTION }\label{s1}
     In his second notebook \cite{4}, S.Ramanujan proved the
following theorem.

\begin{theorem}\label{thm1}
(\cite{4}, p.325;\cite{2}). \slshape Let $\alpha,\beta$ and $\gamma$
denote the roots of the cubic equation\upshape

\begin{equation}\label{1}
x^3-ax^2+bx-1=0.
\end{equation}

\slshape Then, for a suitable determination of roots,\upshape

\begin{equation}\label{2}
\alpha^{\frac 1 3}+\beta^{\frac 1 3}+\gamma^{\frac 1
3}=(a+6+3t)^{\frac 1 3},
\end{equation}

\slshape and \upshape

\begin{equation}\label{3}
(\alpha\beta)^{\frac 1 3}+(\beta\gamma)^{\frac 1
3}+(\gamma\alpha)^{\frac 1 3}=(b+6+3t),
\end{equation}

\slshape where \upshape

\begin{equation}\label{4}
t^3-3(a+b+3)t-(ab+6(a+b)+9)=0.
\end{equation}
\end{theorem}

     A proof of Theorem 1 can be found in paper \cite{2}.
Evidently, the simplest condition for successful application of
Theorem 1 is the condition

\newpage

\begin{equation}\label{5}
a+b+3=0
\end{equation}

that implies

\begin{equation}\label{6}
t=(ab+6(a+b)+9)^{\frac 1 3}=(ab-9)^{\frac 1 3}.
\end{equation}

If, for real nonzero $r$, to consider

\begin{equation}\label{7}
x_1=-r^{\frac 1 3}\alpha,\;\;x_2=-r^{\frac 1 3}\beta,\;\;
x_3=-r^{\frac 1 3}\gamma
\end{equation}

and denote

\begin{equation}\label{8}
ar^{\frac 1 3}=p,\;\;br^{\frac 2 3}=q
\end{equation}

then in the real case by (\ref{1})-(\ref{8}) we obtain the following
result.

\begin{theorem}\label{thm2}

\slshape Let $p, q, r, \in\mathbb{R},\;\;r\neq 0$ such that\upshape

\begin{equation}\label{9}
pr^{\frac 1 3}+ 3r^{\frac 2 3}+q=0,
\end{equation}

\slshape and let the polynomial\upshape

\begin{equation}\label{10}
x^3+px^2+qx+r
\end{equation}

\slshape have real roots $x_1, x_2, x_3$. Then\upshape

\begin{equation}\label{11}
x_1^{\frac 1 3}+x_2^{\frac 1 3}+x_3^{\frac 1 3}=(-p-6r^{\frac 1
3}+3(9r-pq)^{\frac 1 3})^{\frac 1 3},
\end{equation}

\slshape and \upshape

\begin{equation}\label{12}
(x_1x_2)^{\frac 1 3}+(x_2x_3)^{\frac 1 3}+(x_3x_4)^{\frac 1
3}=(q+6r^{\frac 2 3}-3(9r^2-pqr)^{\frac 1 3})^{\frac 1 3},
\end{equation}
\end{theorem}

Notice that Theorem 2 was proved directly in \cite{6}. Notice also,
that (\ref{12}) can be written in the form

\begin{equation}\label{13}
x_1^{-\frac 1 3}+x_2^{-\frac 1 3}+x_3^{-\frac 1 3}=r^{-\frac 1
3}(-q-6r^{\frac 2 3}+3(9r^2-pqr)^{\frac 1 3})^{\frac 1 3}.
\end{equation}

In connection with Theorem 2 we introduce the following definition.
\newpage

\begin{definition}  Let $p, q, r\in \mathbb{R},\;\;r\neq 0$. The
cubic polynomial (\ref{10}) is called a Ramanujan polynomial (RCP)
if it has real roots and the condition (\ref{9}) is satisfied.
\end{definition}

     In  this paper we study various properties of the RCP and
present some new identities.

\section{SOME EXAMPLES}\label{s2}

\bfseries Example 1.\mdseries The polynomial $x^3-3x^2-6x+8$ is an
RCP with roots $1, -2, 4$. Thus, by (\ref{11}) we have

$$
1-2^{\frac 1 3}+4^{\frac 1 3}=(-9+9\cdot 2^{\frac 1 3})^{\frac 1 3}
$$

and therefore

\begin{equation}\label{14}
\left(\frac 1 9\right)^{\frac 1 3}-\left(\frac 2 9\right)^{\frac 1
3}+\left(\frac 4 9\right)^{\frac 1 3}=(2^{\frac 1 3}-1)^{\frac 1 3}.
\end{equation}

     It is a Ramanujan`s original identity (see\cite{5},p.331).

\bfseries Example 2.\mdseries Ramanujan (\cite{4}, p.326,
\cite{1})offers the following identities:

\begin{equation}\label{15}
\left(\cos\frac{2\pi}{7}\right)^{\frac 1 3}+\left(\cos\frac
{4\pi}{7}\right)^{\frac 1 3}+\left(\cos\frac{8\pi}{7}\right)^{\frac
1 3}=\left(\frac{5-3\cdot 7^{\frac 1 3}}{2} \right)^{\frac 1 3}
\end{equation}

\begin{equation}\label{16}
\left(\cos\frac{2\pi}{9}\right)^{\frac 1 3}+\left(\cos\frac
{4\pi}{9}\right)^{\frac 1 3}+\left(\cos\frac{8\pi}{9}\right)^{\frac
1 3}=\left(\frac{3\cdot 9^{\frac 1 3}-6}{2} \right)^{\frac 1 3}
\end{equation}

Notice that (sce, e.g.,\cite{6})

\begin{equation}\label{17}
\left(x-\cos\frac{2\pi}{7}\right)\left(x-\cos\frac
{4\pi}{7}\right)\left(x-\cos\frac{8\pi}{7}\right) =x^3+\frac 1 2
x^2-\frac 1 2 x-\frac 1 8
\end{equation}

\begin{equation}\label{18}
\left(x-\cos\frac{2\pi}{9}\right)\left(x-\cos\frac
{4\pi}{9}\right)\left(x-\cos\frac{8\pi}{9}\right) =x^3-\frac 3 4
x+\frac 1 8
\end{equation}

and both polynomials are $RCP_s$. Thus, by (\ref{11}) we
obtain(\ref{15}) and (\ref{16}).

     Besides, using (\ref{13}), (\ref{17}) and (\ref{18}) we find
\newpage

\begin{equation}\label{19}
\left(\sec\frac{2\pi}{7}\right)^{\frac 1 3}+\left(\sec\frac
{4\pi}{7}\right)^{\frac 1 3}+\left(\sec\frac{8\pi}{7}\right)^{\frac
1 3}=\left(8-6\cdot 7^{\frac 1 3} \right)^{\frac 1 3}
\end{equation}

\begin{equation}\label{20}
\left(\sec\frac{2\pi}{9}\right)^{\frac 1 3}+\left(\sec\frac
{4\pi}{9}\right)^{\frac 1 3}+\left(\sec\frac{8\pi}{9}\right)^{\frac
1 3}=\left(6(9^{\frac 1 3}-1)\right)^{\frac 1 3}
\end{equation}

\bfseries Example 3. \mdseries Quite recently R.Witula and D.Slota
\cite{7} found, in particular, the following decomposition

$$
x^3+105 x^2+588 x-343=
$$
\begin{equation}\label{21}
=(x-2\sin^6\frac{2\pi}{7}\cos\frac{4\pi}{7})(x-2\sin^6\frac{4\pi}{7}\cos\frac{8\pi}{7})
(x-2\sin^6\frac{8\pi}{7}\cos\frac{2\pi}{7})
\end{equation}

This polynomial is an  RCP. Therefore by (\ref{11}) and (\ref{13})
we obtain the following identities (the first of them is presented
in \cite{7})

$$
\sin^2\frac{2\pi}{7}\left(2\cos\frac{4\pi}{7}\right)^{\frac 1
3}+\sin^2\frac{4\pi}{7}\left(2\cos\frac {8\pi}{7}\right)^{\frac 1
3}+
$$
\begin{equation}\label{22}
+\sin^2\frac{8\pi}{7}\left(2\cos\frac{2\pi}{7}\right)^{\frac 1
3}=-\frac 1 4\left(63(1+7^{\frac 1 3})\right)^{\frac 1 3};
\end{equation}

$$
\csc^2\frac{2\pi}{7}\left(2\sec\frac{4\pi}{7}\right)^{\frac 1
3}+\csc^2\frac{4\pi}{7}\left(2\sec\frac {8\pi}{7}\right)^{\frac 1
3}+
$$
\begin{equation}\label{23}
+\csc^2\frac{8\pi}{7}\left(2\sec\frac{2\pi}{7}\right)^{\frac 1
3}=7\left(441(2-7^{\frac 1 3})\right)^{\frac 1 3}
\end{equation}

\section{SOME PROPERTIES OF $ RCP_s $}\label{s3}

\begin{theorem}\label{thm3}

\slshape If $x^3+px^2+qx+r$ is an RCP with roots ${x_1, x_2, x_3}$
then

1) For any $a\in\mathbb{R},\;\;a\neq 0$, the polynomial \upshape

$$x^3+apx^2+a^2qx+a^3r$$
\slshape is also an RCP with roots\upshape\;\; ${ax_1, ax_2, ax_3}$.

2)\slshape The polynomial\upshape

$$
x^3+qx^2+(pr)x+r^2
$$
\slshape is
also an RCP with roots \upshape

\newpage

$$
\left\{\frac{r}{x_1}, \frac{r}{x_2}, \frac{r}{x_3}\right\}.
$$

3)\slshape The numbers\upshape $$\left\{\frac{r^{\frac 2 3}}{x_1},
\frac{r^{\frac 2 3}}{x_2}, \frac{r^{\frac 2 3}}{x_3}\right\}$$
\slshape are a permutation of the numbers \upshape $$\left\{r^{\frac
1 3}-x_1,\;r^{\frac 1 3}-x_2,\;r^{\frac 1 3}-x_3\right\}.$$

4) $\frac{pq}{r}\leq \frac 9 4$.
\end{theorem}

\slshape Proof.\upshape   1)-2)  Straightforward.

3) Let

$$
(x-(x_1-r^{\frac 1 3}))(x-(x_2-r^{\frac 1 3}))(x-(x_3-r^{\frac 1
3}))=x^3+p_1x^2+q_1x+r_1
$$

After some simple calculations we have

$$
p_1=p+3r^{\frac 1 3}
$$
$$
q_1=q+2pr^{\frac 1 3}+3r^{\frac 2 3}
$$
$$
r_1=r+qr^{\frac 1 3}+pr^{\frac 2 3}+r.
$$

Using (\ref{9}) we find

\begin{equation}\label{24}
p_1=-\frac{q}{r^{\frac 1 3}},\; q_1=pr^{\frac 1 3}, \; r_1=-r
\end{equation}

Now we see that (\ref{9}) is satisfied for numbers (\ref{24}):

$$
p_1r_1^{\frac 1 3}+3r_1^{\frac 2 3}+q_1=0.
$$

i.e. the polynomial $x^3+p_1x^2+q_1x+r_1$ is an RCP with the roots
$x_1-r^{\frac 1 3}, \; x_2-r^{\frac 1 3}, \; x_3-r^{\frac 1 3}$.

According to 1) for $a=-r^{\frac 1 3}$ we notice that polynomial
$x^3+qx^2+(pr)x+r^2$ has the roots

$$
\left\{ r^{\frac 1 3}(r^{\frac 1 3}-x_1),\; r^{\frac 1 3}(r^{\frac 1
3}-x_2), \; r^{\frac 1 3}(r^{\frac 1 3}-x_3)\right\}.
$$

On the other hand, by 2) it has the roots $ \left\{\frac
{r}{x_1},\;\frac{r}{x_2},\; \frac{r}{x_3}\right\}$  and 3) follows.

\newpage

4) Notice that, the condition (\ref{9}) yields

\begin{equation}\label{25}
\frac{pq}{r}=-\left(\frac{p}{r^{\frac 1
3}}\right)^2-3\left(\frac{p}{r^{\frac 1 3}}\right)
\end{equation}

Thus,

$$
\frac 9 4 -\frac{pq}{r}=\frac 9 4 +3\frac{p}{r^{\frac 1
3}}+\left(\frac{p}{r^{\frac 1 3}}\right)^2=\left(\frac 3 2+
\frac{p}{r^{\frac 1 3}}\right)^2\geq 0.\blacksquare
$$

\bfseries Example 4. \mdseries By  1) of Theorem 1 together with the
RCP (\ref{17}) we have also the RCP $x^3+x^2-2x-1$ with the roots
$x_1=2\cos\frac{2\pi}{7}, \;x_2=2\cos\frac{4\pi}{7},
\;x_3=2\cos\frac{8\pi}{7}$.

\bfseries Example 5.\mdseries  For the RCP from Example 4  we have
$\frac{r^{\frac 2 3}}{x_1}=\frac 1 2 \sec\frac{2\pi}{7}> 0$ while
$r^{\frac 1 3}-x_1=-1-2\cos\frac{2\pi}{7}< 0$. Therefore, the
permutation in 3) of Theorem 1 is not identical. We verity \slshape
approximately \upshape by a pocket calculator the precise equalities

$$
\frac{1}{x_1}=1+x_3, \; \frac{1}{x_2}=1+x_1, \; \frac{1}{x_3}=1+x_2.
$$

\bfseries Example 6. \mdseries Analogously, by (\ref{18}) we have
the RCP $x^3-3x+1$ with the roots $x_1=2\cos\frac{2\pi}{9}, \;
x_2=2\cos\frac{4\pi}{9}, \; x_3=2\cos\frac{8\pi}{9}$. Here
$\frac{r^{\frac 2 3}}{x_1}=\frac 1 2\csc\frac{2\pi}{9}> 0$ while
$r^{\frac 1 3}-x_1=1-2\cos\frac{2\pi}{9}< 0$. Therefore, the
permutation in

3) of Theorem 1 again is not identical. We verity as above that

$$
\frac{1}{x_1}=1-x_2, \;\frac{1}{x_2}= 1-x_3, \;\frac{1}{x_3}=1-x_1.
$$

\section{A NEW FORMULA AND IDENTITIES}\label{s4}

\begin{theorem}\label{thm4}
\slshape In the conditions of Theorem 3 the following formula
holds\upshape

\begin{equation}\label{26}
\left(\frac{x_1}{x_2}\right)^{\frac 1
3}+\left(\frac{x_2}{x_1}\right)^{\frac 1
3}+\left(\frac{x_1}{x_3}\right)^{\frac 1
3}+\left(\frac{x_3}{x_1}\right)^{\frac 1
3}+\left(\frac{x_2}{x_3}\right)^{\frac 1
3}+\left(\frac{x_3}{x_2}\right)^{\frac 1
3}=\left(\frac{pq}{r}-9\right)^{\frac 1 3}
\end{equation}

\end{theorem}

\slshape  Proof.\upshape Denote X the left hand side of (\ref{26}).
Notice that, in \cite{6} it was shown that (\ref{9}) yields the
following identity for the roots of RCP (\ref{10}):

\begin{equation}\label{27}
X=\left(\frac{x_1}{x_2}+\frac{x_2}{x_1}+\frac{x_1}{x_3}+\frac{x_3}{x_1}
+\frac{x_2}{x_3}+\frac{x_3}{x_2}-6\right)^{\frac 1 3}.
\end{equation}

Thus, we have
\newpage

$$
X=(-\frac 1 r
(x_1^2x_2+x_2^2x_1+x_1^2x_3+x_3^2x_1+x_2^2x_3+x_3^2x_2)-6)^\frac 1
3=
$$
\begin{equation}\label{28}
=(-\frac{1}{3r}((x_1+x_2+x_3)^3-(x_1^3+x_2^3+x_3^3)-6x_1x_2x_3)-6)^{\frac
1 3}.
\end    {equation}

As is well known (see, e.g. \cite{3})

\begin{equation}\label{29}
x_1^3+x_2^3+x_3^3=\sigma_1^3-3\sigma_1\sigma_2+3\sigma_3,
\end{equation}

where $\sigma_i=\sigma_i(x_1, x_2, x_3),\;\;i=1,2,3,$ are the
elementary symmetric polynomials. For RCP (\ref{10}) we have

\begin{equation}\label{30}
\sigma_1=-p,\;\;\sigma_2=q,\;\;\sigma_3=-r.
\end{equation}

Now by (\ref{28})-(\ref{30}) we complete the proof:

$$
X=\left(-\frac{1}{3r}(-p^3+p^3-3pq+3r+6r)-6\right)^{\frac 1
3}=\left(\frac{pq}{r}-9\right)^{\frac 1 3}. \blacksquare
$$

\bfseries Example 7.\mdseries  Using the RCP (\ref{17}),(\ref{18})
by Theorem 4 we find the following identities:

$$
\left(\frac{\cos\frac{2\pi}{7}}{\cos\frac{4\pi}{7}}\right)^{\frac 1
3}+\left(\frac{\cos\frac{4\pi}{7}}{\cos\frac{2\pi}{7}}\right)^{\frac
1
3}+\left(\frac{\cos\frac{2\pi}{7}}{\cos\frac{8\pi}{7}}\right)^{\frac
1
3}+\left(\frac{\cos\frac{8\pi}{7}}{\cos\frac{2\pi}{7}}\right)^{\frac
1 3}+
$$
\begin{equation}\label{31}
+\left(\frac{\cos\frac{4\pi}{7}}{\cos\frac{8\pi}{7}}\right)^{\frac 1
3}+\left(\frac{\cos\frac{8\pi}{7}}{\cos\frac{4\pi}{7}}\right)^{\frac
1 3}=-7^{\frac 1 3};
\end{equation}

$$
\left(\frac{\cos\frac{2\pi}{9}}{\cos\frac{4\pi}{9}}\right)^{\frac 1
3}+\left(\frac{\cos\frac{4\pi}{9}}{\cos\frac{2\pi}{9}}\right)^{\frac
1
3}+\left(\frac{\cos\frac{2\pi}{9}}{\cos\frac{8\pi}{9}}\right)^{\frac
1
3}+\left(\frac{\cos\frac{8\pi}{9}}{\cos\frac{2\pi}{9}}\right)^{\frac
1 3}+
$$
\begin{equation}\label{32}
+\left(\frac{\cos\frac{4\pi}{9}}{\cos\frac{8\pi}{9}}\right)^{\frac 1
3}+\left(\frac{\cos\frac{8\pi}{9}}{\cos\frac{4\pi}{9}}\right)^{\frac
1 3}=-9^{\frac 1 3}.
\end{equation}

\bfseries Example 8.\mdseries  Using RCP(\ref{21}) we find

$$
\left(\frac{\sin\frac{2\pi}{7}}{\sin\frac{4\pi}{7}}\right)^2
\left(\frac{\cos\frac{4\pi}{7}}{\cos\frac{8\pi}{7}}\right)^{\frac 1
3}+\left(\frac{\sin\frac{4\pi}{7}}{\sin\frac{2\pi}{7}}\right)^2
\left(\frac{\cos\frac{8\pi}{7}}{\cos\frac{4\pi}{7}}\right)^{\frac 1
3}+
$$
\newpage

$$
+\left(\frac{\sin\frac{2\pi}{7}}{\sin\frac{8\pi}{7}}\right)^2
\left(\frac{\cos\frac{4\pi}{7}}{\cos\frac{2\pi}{7}}\right)^{\frac 1
3}+\left(\frac{\sin\frac{8\pi}{7}}{\sin\frac{2\pi}{7}}\right)^2
\left(\frac{\cos\frac{2\pi}{7}}{\cos\frac{4\pi}{7}}\right)^{\frac 1
3}+
$$
\begin{equation}\label{33}
+\left(\frac{\sin\frac{4\pi}{7}}{\cos\frac{8\pi}{7}}\right)^2
\left(\frac{\cos\frac{8\pi}{7}}{\cos\frac{2\pi}{7}}\right)^{\frac 1
3}+\left(\frac{\sin\frac{8\pi}{7}}{\sin\frac{4\pi}{7}}\right)^2
\left(\frac{\cos\frac{2\pi}{7}}{\cos\frac{8\pi}{7}}\right)^{\frac 1
3}=-3\cdot 7^{\frac 1 3}.
\end{equation}

\section{ON RCP WITH THE SAME VALUE OF $\frac{pq}{r}$}\label{s5}

\begin{theorem}\label{thm5}
\slshape  If for two $RCP_s$

$$
y^3+p_1y^2+q_1y+r_1, \;\;\;z^3+p_2z^2+q_2z+r_2
$$

the following condition holds

$$
\frac{p_1q_1}{r_1}=\frac{p_2q_2}{r_2}
$$

then for its roots $y_1, y_2, y_3;\;\;z_1, z_2, z_3$ the numbers

$$
\left\{\frac{y_1}{y_2},\;\frac{y_2}{y_1},\;\frac{y_1}{y_3},\;\frac{y_3}{y_1},
\;\frac{y_2}{y_3},\;\frac{y_3}{y_2}\right\}
$$

are a permutation of the numbers

$$
\left\{\frac{z_1}{z_2},\;\frac{z_2}{z_1},\;\frac{z_1}{z_3},\;\frac{z_3}{z_1},
\;\frac{z_2}{z_3},\;\frac{z_3}{z_2}\right\}
$$
\upshape
\end{theorem}

\slshape  Proof.\upshape  Evidently, it is sufficient to prove that
if $x_1, x_2, x_3$ are roots of RCP (\ref{10}) then the numbers

\begin{equation}\label{34}
\xi_1=\frac{x_1}{x_2}+(\frac{x_1}{x_2})^{-1},
\;\;\xi_2=\frac{x_1}{x_3}+(\frac{x_1}{x_3})^{-1},\;\;
\xi_3=\frac{x_2}{x_3}+(\frac{x_2}{x_3})^{-1}
\end{equation}

depend on $\frac{pq}{r}$ only.

For elementary symmetric polynomials of $\xi_1, \xi_2, \xi_3$ we
have
\begin{equation}\label{35}
\xi_1+\xi_2+\xi_3=-\frac 1 r (x_1^2x_2+x_2^2x_1+x_1^2x_3+x_3^2x_1+
x_2^2x_3+x_3^2x_1)
\end{equation}

and as above ((\ref{28})-(\ref{29})) we find
\begin{equation}\label{36}
\xi_1+\xi_2+\xi_3=\frac{pq}{r}-3.
\end{equation}

$$
\xi_1\xi_2+\xi_1\xi_3+\xi_2\xi_3=\frac{1}{r^2}(-r(x_1^3+x_2^3+x_3^3)+
(x_1x_2)^3+(x_1x_3)^3+(x_2x_3)^3-
$$
$$
-r(x_1x_2^2+x_2x_1^2+x_1^2x_3+x_3^2x_1+x_2^2x_3+x_3^2x_2)).
$$

\newpage

Using (\ref{35})-(\ref{36}) and taking into account that by
(\ref{29}),(\ref{30})
$$
x_1^3+x_2^3+x_3^3=-p^3+3pq-3r,
$$
$$
(x_1x_2)^3+(x_1x_3)^3+(x_2x_3)^3=-(p*)^3+3p^*q^*-3r^*,
$$

where

$$
p^*=-(x_1x_2+x_1x_3+x_2x_3)=-q,
$$
$$
q^*=x_1^2x_2x_3+x_1x_2^2x_3+x_1x_2x_3^2=x_1x_2x_3(x_1+x_2+x_3)=pr,
$$
$$
r^*=-(x_1x_2x_3)^2=-r^2,
$$

such that

$$
(x_1x_2)^3+(x_1x_3)^3+(x_2x_3)^3=q^3-3pqr+3r^2,
$$

we have
\begin{equation}\label{37}
\xi_1\xi_2+\xi_1\xi_3+\xi_2\xi_3=4\frac{pq}{r}-24.
\end{equation}

It is interesting that, \slshape only \upshape for proof that
$\xi_1\xi_2\xi_3$ depends on $\frac{pq}{r}$, we use condition
(\ref{9}). Indeed, by (\ref{34}) we have

$$
\xi_1\xi_2\xi_3=\frac{1}{r^2}(x_1^2+x_2^2+x_3^2)((x_1x_2)^2+(x_1x_3)^2+(x_2x_3)^2)-1.
$$

Taking into account that

$$
x_1^2+x_2^2+x_3^2=p^2-2q
$$

$$
(x_1x_2)^2+(x_1x_3)^2+(x_2x_3)^2=(p^*)^2-2q^*=q^2-2pr
$$

we find

\begin{equation}\label{38}
\xi_1\xi_2\xi_3=\left(\frac{pq}{r}\right)^2+4\left(\frac{pq}{r}\right)-1
-2\frac{p^3r+q^3}{r^2}.
\end{equation}

Since by (\ref{9})
$$
pr^{\frac 1 3}+q=-3r^{\frac 2 3}
$$

then we have
\newpage

$$
\frac{p^3r+q^3}{r^2}=\frac{1}{r^2}(pr^{\frac 1 3}+q)(p^2r^{\frac 2
3}-pqr^{\frac 1 3}+q^2)=
$$
$$
=-\frac{3}{r^{\frac 4 3}}(pr^{\frac 1 3}+q)^2-3pqr^{\frac 1 3})=
$$
$$
=-\frac{3}{r^{\frac 4 3}}(9r^{\frac 4 3}-3pqr^{\frac 1 3}).
$$

Now according to (\ref{38}) we find

\begin{equation}\label{39}
\xi_1\xi_2\xi_3=\left(\frac{pq}{r}-7\right)^2+4.
\end{equation}

Thus, by (\ref{36}), (\ref{37}) and (\ref{39}) the numbers
$\xi_1,\xi_2,\xi_3 $ (\ref{34}) are the roots of the following
polynomial

\begin{equation}\label{40}
\xi^3-\left(\frac{pq}{r}-3\right)\xi^2+4\left(\frac{pq}{r}-6\right)\xi
-\left(\left(\frac{pq}{r}-7\right)^2+4\right)
\end{equation}

and, consequently, they depend on $\frac{pq}{r}$ only.$\blacksquare$

Notice that, using (\ref{25}) we see that the numbers
$\frac{x_1}{x_2}, \frac{x_2}{x_1}, \frac{x_1}{x_3}, \frac{x_3}{x_1},
\frac{x_2}{x_3}, \frac{x_3}{x_2}$ depend even on $\frac{p}{r^{\frac
1 3}}$ only.

\bfseries Example 9. \mdseries  Together with decomposition
(\ref{17}) with $\frac{pq}{r}=2$ in \cite{7} the following
decomposition was found:

$$
x^3+7x^2-98x-343=
$$
$$
\left(x-128\cos\frac{2\pi}{7}
\left(\sin\frac{2\pi}{7}\sin\frac{8\pi}{7}\right)^3\right)
\left(x-128\cos\frac{4\pi}{7}
\left(\sin\frac{2\pi}{7}\sin\frac{4\pi}{7}\right)^3\right)\times
$$
$$
\times \left(x-128\cos\frac{8\pi}{7}
\left(\sin\frac{4\pi}{7}\sin\frac{8\pi}{7}\right)^3\right).
$$

It is an RCP with also $\frac{pq}{r}=2$. Hence, by Theorem 5 the
numbers

$$
\left\{\frac{\cos\frac{2\pi}{7}}{\cos\frac{4\pi}{7}}\left(\frac{\sin\frac{8\pi}{7}}
{\sin\frac{4\pi}{7}}\right)^3,\;\;\frac{\cos\frac{4\pi}{7}}{\cos\frac{2\pi}{7}}
\left(\frac{\sin\frac{4\pi}{7}}{\sin\frac{8\pi}{7}}\right)^3,\;\;
\frac{\cos\frac{2\pi}{7}}{\cos\frac{8\pi}{7}}\left(\frac{\sin\frac{2\pi}{7}}
{\sin\frac{4\pi}{7}}\right)^3,\right.
$$

$$
\left.
\frac{\cos\frac{8\pi}{7}}{\cos\frac{2\pi}{7}}\left(\frac{\sin\frac{4\pi}{7}}
{\sin\frac{2\pi}{7}}\right)^3,\;\;\frac{\cos\frac{4\pi}{7}}{\cos\frac{8\pi}{7}}
\left(\frac{\sin\frac{2\pi}{7}}{\sin\frac{8\pi}{7}}\right)^3,\;\;
\frac{\cos\frac{8\pi}{7}}{\cos\frac{4\pi}{7}}\left(\frac{\sin\frac{8\pi}{7}}
{\sin\frac{2\pi}{7}}\right)^3\right\}
$$

are a permutation of the numbers

$$
\left\{\frac{\cos\frac{2\pi}{7}}{\cos\frac{4\pi}{7}},\;\;\;
\frac{\cos\frac{4\pi}{7}}{\cos\frac{2\pi}{7}},\;\;\;
\frac{\cos\frac{2\pi}{7}}{\cos\frac{8\pi}{7}},\;\;\;
\frac{\cos\frac{8\pi}{7}}{\cos\frac{2\pi}{7}},\;\;\;
\frac{\cos\frac{4\pi}{7}}{\cos\frac{8\pi}{7}},\;\;\;
\frac{\cos\frac{8\pi}{7}}{\cos\frac{4\pi}{7}}\right\}.
$$

\newpage

To get the corresponding equalities, it is sufficient to verify
\slshape approximately \upshape by a pocket calculator that

$$
\frac{\cos\frac{2\pi}{7}}{\cos\frac{4\pi}{7}}\left(\frac{\sin\frac{8\pi}{7}}
{\sin\frac{4\pi}{7}}\right)^3=\frac{\cos\frac{4\pi}{7}}{\cos\frac{8\pi}{7}},\;\;
\frac{\cos\frac{2\pi}{7}}{\cos\frac{8\pi}{7}}\left(\frac{\sin\frac{2\pi}{7}}
{\sin\frac{4\pi}{7}}\right)^3=\frac{\cos\frac{4\pi}{7}}{\cos\frac{2\pi}{7}},
$$
\begin{equation}\label{41}
\frac{\cos\frac{4\pi}{7}}{\cos\frac{8\pi}{7}}\left(\frac{\sin\frac{2\pi}{7}}
{\sin\frac{8\pi}{7}}\right)^3=
\frac{\cos\frac{8\pi}{7}}{\cos\frac{2\pi}{7}},
\end{equation}
etc.

Notice that, in \cite{7} R.Witula and D.Slota found some other
$RCP_s$ with roots which depend on sines or cosines of arguments
$\frac{2\pi}{7},\;\;\frac{4\pi}{7},\;\;\frac{8\pi}{7}$. However, it
is interesting that, as we verified, for them $\frac{pq}{r}$ assumes
\slshape only \upshape the values $2, -40, -180$.


\begin{thebibliography}{7}
\bibitem{1} 1. B.C.~Berndt,Ramanujan`s
Notebooks, Part IV,Springer-Verlag, New York,1994.
\bibitem{2} 2. B.C.~Berndt and S.~Bhargava, Ramanujan-for
Lowbraws,  \slshape American Mathematical Monthly \upshape
,\bfseries 100 \mdseries (1993),644-656.
\bibitem{3} 3. A.G.~Kurosh.
\slshape Lectures in General Algebra,\upshape Chelsea,1963.
\bibitem{4} 4. S.~Ramanujan,
\slshape Nitebooks \upshape     (2 volumes), Tata Institute of
Fundamental Research, Bombay, (1957).
\bibitem{5} 5. S.~Ramanujan, \slshape Collected Papers, \upshape
  Chelsea,New York,1962.
\bibitem{6} 6. V.S.~Shevelev. Three Ramanujan`s Formulas,
\slshape Kvant   \upshape \bfseries    6 \mdseries ,(1988), 52-55
(in Russian). Translation into English in AMS:Kvant Selecta,
\slshape   Mathem.World,\upshape vol.\bfseries 14 \mdseries (1999),
139-144.
\bibitem{7} 7. R.~Witula and D.~Slota, New Ramanujan - Type Formulas
and Quasi - Fibonacci Numbers of Order   7, \slshape Integer
Sequences,\upshape  \bfseries 10 \mdseries  (2007), Article 07.5.6.
\end{thebibliography}
\end{document}